\documentclass{article}
\usepackage{amssymb,amsmath,theorem,euscript}

\newcounter{sec}

\def\SS{\smallskip}

\newcounter{punct}[sec]

\def\punct{\refstepcounter{punct}{\arabic{sec}.\arabic{punct}.  }}

\def\COUNTERS{\addtocounter{sec}{1}
              \setcounter{punct}{0}
          \setcounter{equation}{0}
          \setcounter{theorem}{0}
                  }

\newtheorem{theorem}{Theorem}[sec]
\newtheorem{proposition}[theorem]{Proposition}

\begin{document}

 \def\ov{\overline}
\def\wt{\widetilde}
 \newcommand{\rk}{\mathop {\mathrm {rk}}\nolimits}
\newcommand{\Aut}{\mathop {\mathrm {Aut}}\nolimits}
\newcommand{\Out}{\mathop {\mathrm {Out}}\nolimits}
\renewcommand{\Re}{\mathop {\mathrm {Re}}\nolimits}
\def\Br{\mathrm {Br}}

\def\SL{\mathrm {SL}}
\def\SU{\mathrm {SU}}
\def\GL{\mathrm {GL}}
\def\U{\mathrm U}
\def\OO{\mathrm O}
 \def\Sp{\mathrm {Sp}}
 \def\SO{\mathrm {SO}}
\def\SOS{\mathrm {SO}^*}
 \def\Diff{\mathrm{Diff}}
 \def\Vect{\mathfrak{Vect}}
\def\PGL{\mathrm {PGL}}
\def\PU{\mathrm {PU}}
\def\PSL{\mathrm {PSL}}
\def\Symp{\mathrm{Symp}}
\def\End{\mathrm{End}}
\def\Mor{\mathrm{Mor}}
\def\Aut{\mathrm{Aut}}
 \def\PB{\mathrm{PB}}
 \def\cA{\mathcal A}
\def\cB{\mathcal B}
\def\cC{\mathcal C}
\def\cD{\mathcal D}
\def\cE{\mathcal E}
\def\cF{\mathcal F}
\def\cG{\mathcal G}
\def\cH{\mathcal H}
\def\cJ{\mathcal J}
\def\cI{\mathcal I}
\def\cK{\mathcal K}
 \def\cL{\mathcal L}
\def\cM{\mathcal M}
\def\cN{\mathcal N}
 \def\cO{\mathcal O}
\def\cP{\mathcal P}
\def\cQ{\mathcal Q}
\def\cR{\mathcal R}
\def\cS{\mathcal S}
\def\cT{\mathcal T}
\def\cU{\mathcal U}
\def\cV{\mathcal V}
 \def\cW{\mathcal W}
\def\cX{\mathcal X}
 \def\cY{\mathcal Y}
 \def\cZ{\mathcal Z}
\def\0{{\ov 0}}
 \def\1{{\ov 1}}
 \def\frA{\mathfrak A}
 \def\frB{\mathfrak B}
\def\frC{\mathfrak C}
\def\frD{\mathfrak D}
\def\frE{\mathfrak E}
\def\frF{\mathfrak F}
\def\frG{\mathfrak G}
\def\frH{\mathfrak H}
\def\frI{\mathfrak I}
 \def\frJ{\mathfrak J}
 \def\frK{\mathfrak K}
 \def\frL{\mathfrak L}
\def\frM{\mathfrak M}
 \def\frN{\mathfrak N} \def\frO{\mathfrak O} \def\frP{\mathfrak P} \def\frQ{\mathfrak Q} \def\frR{\mathfrak R}
 \def\frS{\mathfrak S} \def\frT{\mathfrak T} \def\frU{\mathfrak U} \def\frV{\mathfrak V} \def\frW{\mathfrak W}
 \def\frX{\mathfrak X} \def\frY{\mathfrak Y} \def\frZ{\mathfrak Z} \def\fra{\mathfrak a} \def\frb{\mathfrak b}
 \def\frc{\mathfrak c} \def\frd{\mathfrak d} \def\fre{\mathfrak e} \def\frf{\mathfrak f} \def\frg{\mathfrak g}
 \def\frh{\mathfrak h} \def\fri{\mathfrak i} \def\frj{\mathfrak j} \def\frk{\mathfrak k} \def\frl{\mathfrak l}
 \def\frm{\mathfrak m} \def\frn{\mathfrak n} \def\fro{\mathfrak o} \def\frp{\mathfrak p} \def\frq{\mathfrak q}
 \def\frr{\mathfrak r} \def\frs{\mathfrak s} \def\frt{\mathfrak t} \def\fru{\mathfrak u} \def\frv{\mathfrak v}
 \def\frw{\mathfrak w} \def\frx{\mathfrak x} \def\fry{\mathfrak y} \def\frz{\mathfrak z} \def\frsp{\mathfrak{sp}}
 \def\bfa{\mathbf a} \def\bfb{\mathbf b} \def\bfc{\mathbf c} \def\bfd{\mathbf d} \def\bfe{\mathbf e} \def\bff{\mathbf f}
 \def\bfg{\mathbf g} \def\bfh{\mathbf h} \def\bfi{\mathbf i} \def\bfj{\mathbf j} \def\bfk{\mathbf k} \def\bfl{\mathbf l}
 \def\bfm{\mathbf m} \def\bfn{\mathbf n} \def\bfo{\mathbf o} \def\bfp{\mathbf p} \def\bfq{\mathbf q} \def\bfr{\mathbf r}
 \def\bfs{\mathbf s} \def\bft{\mathbf t} \def\bfu{\mathbf u} \def\bfv{\mathbf v} \def\bfw{\mathbf w} \def\bfx{\mathbf x}
 \def\bfy{\mathbf y} \def\bfz{\mathbf z} \def\bfA{\mathbf A} \def\bfB{\mathbf B} \def\bfC{\mathbf C} \def\bfD{\mathbf D}
 \def\bfE{\mathbf E} \def\bfF{\mathbf F} \def\bfG{\mathbf G} \def\bfH{\mathbf H} \def\bfI{\mathbf I} \def\bfJ{\mathbf J}
 \def\bfK{\mathbf K} \def\bfL{\mathbf L} \def\bfM{\mathbf M} \def\bfN{\mathbf N} \def\bfO{\mathbf O} \def\bfP{\mathbf P}
 \def\bfQ{\mathbf Q} \def\bfR{\mathbf R} \def\bfS{\mathbf S} \def\bfT{\mathbf T} \def\bfU{\mathbf U} \def\bfV{\mathbf V}
 \def\bfW{\mathbf W} \def\bfX{\mathbf X} \def\bfY{\mathbf Y} \def\bfZ{\mathbf Z} \def\bfw{\mathbf w}
 \def\R {{\mathbb R }} \def\C {{\mathbb C }} \def\Z{{\mathbb Z}} \def\H{{\mathbb H}} \def\K{{\mathbb K}}
 \def\N{{\mathbb N}} \def\Q{{\mathbb Q}} \def\A{{\mathbb A}} \def\T{\mathbb T} \def\P{\mathbb P} \def\G{\mathbb G}
 \def\bbA{\mathbb A} \def\bbB{\mathbb B} \def\bbD{\mathbb D} \def\bbE{\mathbb E} \def\bbF{\mathbb F} \def\bbG{\mathbb G}
 \def\bbI{\mathbb I} \def\bbJ{\mathbb J} \def\bbL{\mathbb L} \def\bbM{\mathbb M} \def\bbN{\mathbb N} \def\bbO{\mathbb O}
 \def\bbP{\mathbb P} \def\bbQ{\mathbb Q} \def\bbS{\mathbb S} \def\bbT{\mathbb T} \def\bbU{\mathbb U} \def\bbV{\mathbb V}
 \def\bbW{\mathbb W} \def\bbX{\mathbb X} \def\bbY{\mathbb Y} \def\kappa{\varkappa} \def\epsilon{\varepsilon}
 \def\phi{\varphi} \def\le{\leqslant} \def\ge{\geqslant}

\def\UU{\bbU}
\def\Mat{\mathrm{Mat}}
\def\tto{\rightrightarrows}

\def\Gr{\mathrm{Gr}}

\def\graph{\mathrm{graph}}

\def\la{\langle}
\def\ra{\rangle}

\begin{center}
{\bf\Large Multi-operator colligations
 and

multivariate characteristic functions}

\bigskip

{\sc\large Yury A. Neretin%
\footnote{Supported by grants FWF, P22122 and P19064.}}

\end{center}

{\small In the spectral theory of non-self-adjoint
operators there is  a well-known operation 
of product of operator colligations. Many similar operations appear in the
 theory of infinite-dimensional groups as multiplications
of double cosets. We
construct characteristic functions for such double cosets and
 get semigroups of matrix-valued functions in matrix balls.}

\section{Introduction}

\COUNTERS

{\bf \punct  Operator colligations.}
Originally, operator colligations
and characteristic functions appeared in the spectral theory of
non-selfadjoint operators in 1946--55 in works of M.~S.~Livshits
and V.~P.~Potapov (see, \cite{Liv1}, \cite{Liv2},
\cite{Pot}, for expositions, see \cite{Goh}, \cite{AD}).
 We discuss  definitions
in a minimal generality and do not touch spectral theory 
and related function theory (see \cite{Nik}).

We say that a unitary operator $U$ in a Hilbert space is
{\it finite}
if rank of $U-1$ is finite.

Consider a finite dimensional Euclidean
 space $\cH$, and an infinite-dimensional Hilbert space $\cK\simeq \ell_2$.
 An {\it operator colligation} is a finite unitary operator
 \begin{equation}
 \frA=\begin{pmatrix}A&B\\C&D\end{pmatrix}: \cH\oplus \cK\to\cH\oplus \cK 
\end{equation}
determined up to a conjugation
\begin{equation}
  \begin{pmatrix}A&B\\C&D\end{pmatrix} \sim
 \begin{pmatrix}1&0\\0&U\end{pmatrix} 
\begin{pmatrix}A&B\\C&D\end{pmatrix} \begin{pmatrix}1&0\\0&U^{-1}\end{pmatrix}
\label{eq:colligation-1}
, \end{equation}
where $U:\cK\to\cK$ is a unitary operator.

 A {\it product of operator colligations}
$$ \frA= \begin{pmatrix}A&B\\C&D\end{pmatrix},\qquad
 \frP=\begin{pmatrix}P&Q\\R&T\end{pmatrix} $$
is given by the formula
\begin{equation} 
\begin{pmatrix}A&B\\C&D\end{pmatrix}
 \circ \begin{pmatrix}P&Q\\R&T\end{pmatrix}
 := \begin{pmatrix}A&B&0\\C&D&0\\0&0&1\end{pmatrix}
 \begin{pmatrix}P&0&Q\\0&1&0\\R&0&T\end{pmatrix} =
\begin{pmatrix} AP&B&AQ\\ CP&D&CQ\\ R&0&T \end{pmatrix}
.
\label{eq:product}
\end{equation}

We get an operator
 $\cH\oplus \cK\oplus\cK\to\cH\oplus \cK\oplus\cK$.
 But all infinite-dimensional
 separable Hilbert spaces are isomorphic,
  we identify $\cK\oplus\cK\simeq \cK$
in  arbitrary way and come to an operator
$\cH\oplus \cK\to\cH\oplus \cK$.

\begin{theorem}
\label{th:1-1}
 a) The multiplication
 $\circ$ is well-defined on the set of operator
 colligations.

 \smallskip

 b) The operation $\circ$ is associative.
\end{theorem}

 Verification is straightforward.

\smallskip


 {\bf\punct Characteristic functions.}
Denote by $\ov\C=\C\cup \infty$ the Riemann sphere.
We define a {\it characteristic function} of an operator colligation
$\frA= \begin{pmatrix}A&B\\C&D\end{pmatrix}$ by
$$
\chi(\frA;z)=A+z B(1-zD)^{-1} C, \qquad z\in\ov\C.
$$

\begin{theorem}
\label{th:1-2}
a) $\chi(\frA;z)$ is a rational matrix-valued function.

\smallskip

b) $\chi(\frA\circ\frP;z)=\chi(\frA)\chi(\frP;z)$
\end{theorem}

\begin{theorem}
\label{th:1-3}
a) $\|\chi(\frA;z)\|\le 1$ for $|z|<1$.

\smallskip

b) $\chi(\frA;z)$ is unitary for $|z|=1$.

\smallskip

c) $\chi(\frA,\ov z^{-1})=\chi(\frA;z)^{*-1}$.
\end{theorem}

Generally speaking, an operator colligation can not be
uniquely reconstructed
from its characteristic function. The reason is the following. The operator colligations
$$
\frA:=\begin{pmatrix}A&B\\C&D \end{pmatrix}, \qquad
\frA:=
\begin{pmatrix}
A&B'&0
\\
C'&D'&0
\\
0&0&L
\end{pmatrix},\qquad\text{where $L$ is unitary},
$$
have the same characteristic function.

For a colligation $\frA:=\begin{pmatrix}A&B\\C&D\end{pmatrix}$
we denote by $\Xi(\frA)$ the set of eigenvalues of $D$ lying on the unit circle
taking in account multiplicities (the multiplicity of $z=1$ is $\infty$).

\smallskip

{\sc Remark.} The spectrum of the block $D$ is contained in the circle $|\lambda|\le1$.
It can be shown that for $|\lambda|<1$, the point $z=\lambda^{-1}$ 
is a pole of the characteristic function. \hfill $\square$

\begin{theorem}
\label{th:1-4}
  Any operator colligation can be uniquely reconstructed from
$\chi(\frA;z)$ and $\Xi(\frA)$.
\end{theorem}

\begin{proposition}
$\Xi(\frA\circ\frP)$ is $\Xi(\frA)\cup\Xi(\frP)$ taking in account multiplicities. 
\end{proposition}

 \smallskip

Denote by $G//L$ conjugacy classes of a group $G$ with respect to
a subgroup $L$.
We can regard $\circ$ as an operations on conjugacy
 classes
 \begin{equation}
\U_{\alpha+n}//\U_n\,\, \times\,\,
\U_{\alpha+m}//\U_m \,\to\,
\U_{\alpha+n+m}//\U_{n+m}
\label{eq:UUU}
 \end{equation}

We also can reject the unitarity condition. Then we come
to a multiplication
 \begin{equation}
\GL_{\alpha+n}//\GL_n\,\, \times\,\,
\GL_{\alpha+m}//\GL_m \,\to\,
\GL_{\alpha+n+m}//\GL_{n+m}
\label{eq:GLGL}
\end{equation}

Then Theorems \ref{th:1-1}, \ref{th:1-2}
survive, Theorem \ref{th:1-3} disappear. Theorem \ref{th:1-4}
 exists in a weaker form  (see, e.g.,\cite{Dym}).

In fact there is lot of operations of this type. An independent
origin is explained in the following two subsections.

Below we prefer to discuss unitary groups and not  $\GL$.


\smallskip

{\bf\punct Another stand-point. Infinite-dimensional classical groups,}
see \cite{OlshGB}, \cite{Ner-book}.
Consider some series $G(n)/K(n)$ of Riemannian symmetric spaces, 
say $\U(n)/\OO(n)$. Here $\U(n)$ is the unitary group and $\OO(n)$
is the real orthogonal group.
Consider double cosets%
\footnote{Let $G$ be a group, $K$ a subgroup. A {\it double coset}
is a set of the form $KgK\subset G$. The notation for sets of all double
cosets is $K\setminus G/K$.}
$$
K(n-\alpha)\setminus G(n)/K(n-\alpha)=
\OO(n-\alpha)\setminus\U(n)/\OO(n-\alpha)
.
$$
They form a 'hypergroup' in the following sense. Let $g\in\U(n)$.
Denote by $\mu_g$ the natural probability measure
on the set  $K(n-\alpha)g K(n-\alpha)$. Consider the convolution
of measures
$$
\mu_{g_1} * \mu_{g_2}=\int \mu_h\,d\lambda(h)
,$$  
where $\lambda(h)$ is a probability measure on  
$K(n-\alpha)\setminus G/K(n-\alpha)$. Thus we get a map
\begin{multline*}
K(n-\alpha)\setminus G(n)/K(n-\alpha) 
\times K(n-\alpha)\setminus G(n)/K(n-\alpha)
\to\\ \to
\Bigl\{ \text{probability measures on
$K(n-\alpha)\setminus G(n)/K(n-\alpha)$
}
\Bigr\}
\end{multline*}
Explicit description of this hypergroup even for
$\alpha=0$ is complicated%
\footnote{As far as I know formulas for 
$K(n)\setminus G(n)/K(n)$  exist only for rank 1 groups
and for complex groups.
See formulas for 
$\OO(2)\setminus\SL(2,\R)/\OO(2)$
in \cite{Koo}, for $\U(n)\setminus\GL(n,\C)/\U(n)$
in \cite{BG}.}. 

However, if we pass to a limit as $n\to\infty$ (and keep $\alpha$ fixed),
then the measure $\nu$ is concentrated near a single double coset%
\footnote{Apparently, the phenomenon of concentration was
firstly observed in \cite{Olsh-tree}.} $\nu$
and a well-defined operation
\begin{multline}
K(\infty-\alpha)\setminus G(\infty)/K(\infty-\alpha) 
\,\,\times\,\,
K(\infty-\alpha)\setminus G(\infty)/K(\infty-\alpha) 
\to\\ \to
K(\infty-\alpha)\setminus G(\infty)/K(\infty-\alpha) 
\label{eq:dcm}
\end{multline}
Note, that  the last operation  is not a convolution of measures
(because there is no a natural measure on 
$K(\infty-\alpha)gK(\infty-\alpha)$).

For instance, consider the case 
$\OO(\infty-\alpha) \setminus\U(\infty)/\OO(\infty-\alpha)$, 
i.e. finite $(\alpha+\infty)\times(\alpha+\infty)$ unitary
matrices defined up to the equivalence
$$
  \begin{pmatrix}A&B\\C&D\end{pmatrix} \sim
 \begin{pmatrix}1&0\\0&U\end{pmatrix} 
\begin{pmatrix}A&B\\C&D\end{pmatrix} \begin{pmatrix}1&0\\0&V\end{pmatrix}
,
$$
where $U$, $V$ are orthogonal matrices.
The multiplication of double cosets is given by the  formula
(\ref{eq:product}).

\begin{theorem}
(Multiplicativity theorem)
Let $\rho$ be a unitary irreducible representation of
$\U(\infty)$ in a Hilbert space $H$. Denote by $H(\alpha)$ the space
of $\OO(\infty-\alpha)$-fixed vectors in $H$.
Denote by $P(\alpha)$ the projection to $H(\alpha)$. Assume that
at least one subspace $H(\beta)$ is non-zero.
For $g\in\U(n)$ define the operator
$$
\ov\rho_\alpha(g)=P(\alpha)\rho(g):\,\,H(\alpha)\to H(\alpha)
.
$$
Then $\ov\rho_\alpha(g)$ is a function of  a double coset $\frA\ni g$
and
\begin{equation}
\ov\rho_\alpha(\frA_1) \ov\rho_\alpha(\frA_2) =
\ov\rho_\alpha(\frA_1\circ \frA_2) 
\label{eq:multiplicativity}
\end{equation}
\end{theorem}

Note that neither the product (\ref{eq:dcm}), nor multiplicativity theorem
 (\ref{eq:multiplicativity})
have  finite-dimensional analogs.

A description of this multiplication in the terms of characteristic functions
and some additional data is given in \cite{Ner-book}, IX.4.
Olshansky paper \cite{OlshGB} provides us a zoo of such constructions 
related to infinite-dimensional symmetric spaces.


\smallskip

{\bf\punct Purposes of the paper.} 
In \cite{Ner-book} it was observed that multiplications
of double cosets and multiplicativity
theorems exist under rather weak restrictions.
Let $G$ be an infinite-dimensional classical group,
$K$ be its subgroup isomorphic to a complete unitary  group $\U(\infty)$
(or $\OO(\infty)$, $\Sp(\infty)$), the subgroup
is equipped with a weak operator topology.
Then usually there is a multiplication
$$
K\setminus G/K\,\,\times\,\, K\setminus G/K \to K\setminus G/K
$$
and usually the multiplicativity theorem holds.

This produces numerous operations of the type
(\ref{eq:product}). Oue purpose is to transfer such operations to
multiplications of meromorphic multi-variate matrix-valued functions.
We explain the technology in Sections 2-3, it can be applied
 in numerous situations. Some further examples
are discussed  Sections 3-5.

Note that nontrivial representation-theoretical constructions
related to infinite-dimensional non-symmetric pairs $G\supset K$ were considered in
\cite{Ness1}, \cite{Ness2} and \cite{Ner-symm}.

\smallskip

{\bf\punct Notation.} Below 

--- $\Mat(n)$ is the space of $n\times n$ matrices;

--- $A^t$, $A^*$ are transposed matrix and adjoint matrix;

---  $\langle\cdot,\cdot\rangle$, $(\cdot,\cdot)$ are
 the standard inner product  and the standard bilinear form on
 $\C^k$,
$$
\la p,q\ra=\sum_{j=1}^k p_j\ov q_j;\qquad
( p,q)=\sum_{j=1}^k p_j q_j;
$$

--- $\U(n)$, $\U(p,q)$, $\Sp(2n,\C)$, are the usual 
notation for classical groups; $\U(\infty)$ denotes the group of
{\it finite} unitary matrices. 

\smallskip

{\bf Acknowledgements.} I am grateful to S.~L.~Tregub and A~.A.~Rosly
for discussion of this topic.

\section{Multiple colligations}

\COUNTERS

{\bf\punct Multiple colligations.}
We say that an {\it $n$-colligation} is a collection
$\frA$ of unitary $(\alpha+\infty)\times(\alpha+\infty)$-matrices 
$g_j=\begin{pmatrix}a_j&b_j\\c_j&d_j\end{pmatrix}$, where $j=1$, \dots, $n$,
defined up to simultaneous conjugation
\begin{multline*}
\left\{\begin{pmatrix}a_1&b_1\\c_1&d_1\end{pmatrix},\,\dots,\,
\begin{pmatrix}a_n&b_n\\c_n&d_n\end{pmatrix}\right\}
\sim\\ \sim
\left\{\begin{pmatrix} 1&0\\0& u\end{pmatrix}
  \begin{pmatrix}a_1&b_1\\c_1&d_1\end{pmatrix}\begin{pmatrix} 1&0\\0& u\end{pmatrix}^{-1},\,\dots,\,
\begin{pmatrix} 1&0\\0& u\end{pmatrix}
\begin{pmatrix}a_n&b_n\\c_n&d_n\end{pmatrix}
\begin{pmatrix} 1&0\\0& u\end{pmatrix}^{-1}
\right\}
,
\end{multline*}
where $u$ is unitary.

\smallskip


{\bf\punct Product of multiple colligations} is determined element-wise,
$$
\Biggl\{
\begin{pmatrix}
a_j&b_j\\c_j&d_j 
\end{pmatrix}
\Biggr\}
\,\circ\,
\Biggl\{
\begin{pmatrix}
\wt a_j&\wt b_j\\ \wt c_j&\wt d_j 
\end{pmatrix}
\Biggr\}=
\Biggl\{
\begin{pmatrix}
a_j&b_j\\c_j&d_j 
\end{pmatrix}
\,\circ\,
\begin{pmatrix}
\wt a_j&\wt b_j\\ \wt c_j&\wt d_j 
\end{pmatrix}
\Biggr\}
.$$

{\bf\punct Characteristic functions.
\label{ss:char-1}}
For definiteness, set $n=3$. 
Fix a $3\times 3$-matrix $S=\{s_{ij}\}$.
We write the equation
\begin{equation}
\begin{pmatrix}
q_1\\s_{11}x_1+s_{12}x_2+s_{13}x_3  
\\q_2\\ s_{21}x_1+s_{22}x_2+s_{23}x_3
 \\q_3\\ s_{31}x_1+s_{32}x_2+s_{33}x_3
 \end{pmatrix}
=
\begin{pmatrix}
a_1&b_1&0&0 &0&0\\
c_1&d_1&0&0 &0&0 \\
0&0&a_2&b_2 &0&0 \\
0&0&c_2&d_2 &0&0\\
0&0&0&0&a_3&b_3\\
0&0&0&0&c_3&d_3
\end{pmatrix}
\begin{pmatrix}
p_1\\ x_1\\ p_2\\ x_2\\p_2\\x_3
\end{pmatrix}
\label{eq:def-2}
\end{equation}
or
\begin{align}
q_1=a_1p_1+b_1 x_1
\label{eq:1}
\\
s_{11}x_1+s_{12}x_2+s_{13}x_3 = c_1p_1+d_1x_1
\label{eq:2}
\\
q_2=a_2p_2+b_2 x_2
\label{eq:3}
\\
s_{21}x_1+s_{22}x_2+s_{23}x_3 = c_2p_2+d_2x_2
\label{eq:4}
\\
q_3=a_3p_3+b_3 x_3
\label{eq:5}
\\
s_{31}x_1+s_{32}x_2+s_{33}x_3 = c_3p_3+d_3x_3
\label{eq:6}
\end{align}

Next, we exclude 'indeterminantes' $x_1$, $x_2$, $x_3$
from (\ref{eq:2}),  (\ref{eq:4}),  (\ref{eq:6}).
Substituting $x_1$, $x_2$, $x_3$ to
(\ref{eq:1}),  (\ref{eq:3}),  (\ref{eq:5})
  we get
a certain dependence of the form
$$
\begin{pmatrix}
 q_1\\q_2\\q_3
\end{pmatrix}
=
\chi(\frA;S)
\begin{pmatrix}
 p_1\\p_2\\p_3
\end{pmatrix}
,$$
where $\chi(\frA;S)\in\Mat(n\alpha)$ (above $n=3$).

Consider the {\it eigensurface} $\Xi[\frA]$
of $d_1$, $d_2$, $d_3$ in the space of $3\times 3$ matrices 
determined by the equation
\begin{equation}
\det
\begin{pmatrix}
 s_{11}-d_1& s_{12}&s_{13}
\\
s_{21}&s_{22}-d_2&s_{23}
\\
s_{31}&s_{32}&s_{33}-d_3
\end{pmatrix}=0
.
\label{eq:eigensurface}
\end{equation}
(here $d_1$, \dots, $d_3$ are given and $s_{ij}$ are indeterminances%
\footnote{There are two points of view to spectral data of several matrices
$A_j$.
The first one is related to determinantal hypersurfaces, see a survey
of Beauville \cite{Bea}. The second is related to
spectral curves, which are widely
 explored in the theory
of integrable systems, see an introduction of Hitchin \cite{Hit}.
 It seems that characteristic functions give the third point of view.}.
If $S\notin\Xi(\frA)$, then
the equations (\ref{eq:2}),  (\ref{eq:4}),  (\ref{eq:6})
have a unique solution.

Now let $n$ arbitrary.

\begin{theorem}
\label{th:main-1}
a) $\chi(\frA;S)$ is a meromorphic matrix-valued function
on the space of $n\times n$-matrices, whose singularities are contained
in the eigensurface $\Xi[\frA]$.

\smallskip

b) $\chi(\frA;S)$ depends only on the operator colligations but not
on matrices $\begin{pmatrix}a_j&b_j\\c_j&d_j\end{pmatrix}$
themselves.

 \smallskip

c) The following identity holds 
$$ 
\chi(\frA;S)\,\chi(\frP;S)=\chi(\frA\circ\frP;S)
$$
pointwise.

d) If $\|S\|<1$, then $\chi(\frA;S)$ is expanding, i.e., it satisfies
$$
\|\chi(\frA;S)^{-1}\|< 1
$$
If $\|S\|= 1$, then we have $\|\chi(\frA;S)^{-1}\|= 1$.

e) In particular, the matrix valued  function $\chi(\frA;S)^{-1}$
is holomorphic in the matrix ball $\|S\|<1$.

\smallskip

f) If $S$ is unitary, then $\chi(\frA;S)$ is unitary.

\smallskip

g) The following Riemann--Schwarz type identity holds
$$
\chi(\frA;S^{*-1})=\chi(\frA;S)^{\,*-1}
$$

h) Let $\lambda_1$, \dots, $\lambda_n\in\C^*$. Let
$\lambda$ be the diagonal matrix with entries $\lambda_j$, and
$\Lambda$ be $(\alpha+\dots+\alpha)\times (\alpha+\dots+\alpha)$
block diagonal matrix with blocks  $\lambda_j$.
Then
$$
\chi(\frA;\lambda S \lambda^{-1})=
\Lambda \,\chi(\frA;S)\,\Lambda^{-1}
.
$$
\end{theorem}


{\bf\punct Proof of Theorem \ref{th:main-1}.a
\label{ss:th-a}}
If $S$ is not on the eigensurface, then
 the equations (\ref{eq:2}), (\ref{eq:4}), (\ref{eq:6})
have a unique solution $x_1$, $x_2$, $x_3$ for given
$p_1$, $p_2$, $p_3$. Therefore the equations
 (\ref{eq:1}), (\ref{eq:3}), (\ref{eq:5}) uniquely determine $q_1$, $q_2$, $q_3$
from $p_1$, $p_2$, $p_3$.

\smallskip


{\bf\punct Proof of Theorem \ref{th:main-1}.b.}
To be definite, set $n=2$. Consider a multiple colligation
 equivalent to a given one. We write the equation
$$
\begin{pmatrix}
 q_1\\s_{11}x_1+s_{12}x_2\\q_2\\s_{21}x_1+s_{22}x_2
\end{pmatrix}
=
\begin{pmatrix}
 1&0&0&0\\
0&u^{-1}&0&0\\
0&0&1&0\\
0&0&0&u^{-1}
\end{pmatrix}
\begin{pmatrix}
a_1&b_1&0&0\\
c_1&d_1&0&0\\
0&0&a_2&b_2\\
0&0&c_2&d_2 
\end{pmatrix}
\begin{pmatrix}
 1&0&0&0\\
0&u&0&0\\
0&0&1&0\\
0&0&0&u
\end{pmatrix}
\begin{pmatrix}
 p_1\\x_1\\p_2\\x_2
\end{pmatrix}
$$
or
$$
\begin{pmatrix}
 q_1\\s_{11}ux_1+s_{12}ux_2\\q_2\\s_{21}ux_1+s_{22}ux_2
\end{pmatrix}
=
\begin{pmatrix}
a_1&b_1&0&0\\
c_1&d_1&0&0\\
0&0&a_2&b_2\\
0&0&c_2&d_2 
\end{pmatrix}
\begin{pmatrix}
 p_1\\ux_1\\p_2\\ux_2
\end{pmatrix}
$$
We denote
$$ 
y_1=ux_1\qquad y_2=ux_2
$$
and come to the system determining $\chi(\frA;S)$.

\smallskip

\
{\bf\punct Proof of Theorem \ref{th:main-1}.c.
\label{ss:th-c}}
To be definite (and to have finite size of matrices), take $n=2$.
We have
$$
\begin{pmatrix}
 q_1\\s_{11}x_1+s_{12}x_2\\q_2\\s_{21}x_1+s_{22}x_2
\end{pmatrix}
=
\begin{pmatrix}
a_1&b_1&0&0\\
c_1&d_1&0&0\\
0&0&a_2&b_2\\
0&0&c_2&d_2 
\end{pmatrix}
\begin{pmatrix}
 p_1\\x_1\\p_2\\x_2
\end{pmatrix}
$$
$$
\begin{pmatrix}
 p_1\\s_{11}y_1+s_{12}y_2\\p_2\\s_{21}y_1+s_{22}xy_2
\end{pmatrix}
=
\begin{pmatrix}
\wt a_1&\wt b_1&0&0\\
\wt c_1&\wt d_1&0&0\\
0&0&\wt a_2&\wt b_2\\
0&0&\wt c_2&\wt d_2 
\end{pmatrix}
\begin{pmatrix}
 r_1\\y_1\\r_2\\y_2
\end{pmatrix}
$$
Then
$$
\begin{pmatrix}
 q_1\\s_{11}x_1+s_{12}x_2\\s_{11}y_1+s_{12}y_2 \\q_2\\s_{21}x_1+s_{22}x_2\\
s_{21}y_1+s_{22}y_2
\end{pmatrix}
=
\begin{pmatrix}
a_1&b_1&0&0&0&0\\
c_1&d_1&0&0&0&0\\
0&0&1&0&0&0\\
0&0&0& a_2&b_2&0\\
0&0&0& c_2&d_2&0\\ 
0&0&0& 0&0&1
\end{pmatrix}
\begin{pmatrix}
 p_1\\x_1\\ s_{11}y_1+s_{12}y_2   \\p_2\\x_2\\ s_{21}y_1+s_{22}y_2
\end{pmatrix}=
$$
$$
=
\begin{pmatrix}
a_1&b_1&0&0&0&0\\
c_1&d_1&0&0&0&0\\
0&0&1&0&0&0\\
0&0&0& a_2&b_2&0\\
0&0&0& c_2&d_2&0\\ 
0&0&0& 0&0&1
\end{pmatrix}
\begin{pmatrix}
\wt a_1&0&\wt b_1&0&0&0\\
0&1&0&0&0&0\\
\wt c_1&0&\wt d_1&0&0&0\\
0&0&0& \wt a_2&0&\wt b_2\\
0&0&0& 0&1&0\\
0&0&0& \wt c_2&0&\wt d_2 
\end{pmatrix}
\begin{pmatrix}
 r_1\\x_1\\y_1\\r_2\\x_2\\y_2
\end{pmatrix}
$$
We get a product of colligations in the right-hand side.

\SS


{\bf\punct Proof of Theorem \ref{th:main-1}.d-f.
\label{ss:th-d}}
Let us prove d).
Since the big matrix in (\ref{eq:def-2})
is unitary, we have
$$
\|q_1\|^2+\|s_{11}x_1+s_{12}x_2\|^2 +\|q_2\|^2+\|s_{21}x_1+s_{22}x_2\|^2
=\|p_1\|^2+\|x_1\|^2+\|p_2\|^2+\|x_2\|^2
$$
Since
$\|S\|\le 1$, we have
\begin{equation}
\|s_{11}x_1+s_{12}x_2\|^2 +\|s_{21}x_1+s_{22}x_2\|^2
\le\|x_1\|^2+\|x_2\|^2
\label{eq:unitary-1}
\end{equation}
and therefore 
\begin{equation}
\|q_1\|^2+\|q_2\|^2\ge \|p_1\|^2+\|p_2\|^2
\label{eq:unitary-2}
\end{equation}
If $S$ is unitary, then we have $=$ in (\ref{eq:unitary-1}) 
and therefore $=$ in (\ref{eq:unitary-2}). 


\smallskip

{\bf\punct Proof of Theorem \ref{th:main-1}.g.\label{ss:th-g}}
Let $n=2$. Since matrices $\begin{pmatrix}
                            a_1&b_1\\ c_1&d_1
                           \end{pmatrix}$,
$\begin{pmatrix}
                            a_2&b_2\\ c_2&d_2
                           \end{pmatrix}$
are unitary, we can write
(\ref{eq:def-2}) as
$$
\begin{pmatrix}
 q_1\\s_{11}x_1+s_{12}x_2\\q_2\\ s_{21}x_1+s_{22}x_2
\end{pmatrix}=
\begin{pmatrix}
\begin{pmatrix}
a_1^*&c_1^*\\b_1^*&d_1^*
\end{pmatrix}^{-1} 
&
\begin{matrix}
 0&0\\0&0
\end{matrix}
\\
\begin{matrix}
 0&0\\0&0
\end{matrix}
&
\begin{pmatrix}
a_2^*&c_2^*\\b_2^*&d_2^*
\end{pmatrix}^{-1} 
\end{pmatrix}
\begin{pmatrix}
p_1&\\x_1&\\p_2\\x_2
\end{pmatrix}
$$
or
$$
\begin{pmatrix}
p_1&\\x_1&\\p_2\\x_2
\end{pmatrix}
=
\begin{pmatrix}
a_1^*&c_1^*&0&0
\\b_1^*&d_1^*&0&0
\\ 
0&0&a_2^*&c_2^*\\
0&0&b_2^*&d_2^* 
\end{pmatrix}
\begin{pmatrix}
 q_1\\s_{11}x_1+s_{12}x_2\\q_2\\ s_{121}x_1+s_{22}x_2
\end{pmatrix}
$$
Now we change variables
$$
\begin{cases}
y_1=s_{11}x_1+s_{12}x_2
\\
y_2=s_{21}x_1+s_{22}x_2
\end{cases}
\qquad
\begin{cases}
x_1=\sigma_{11}y_1+\sigma_{12}y_2
\\
x_2=\sigma_{21}y_1+\sigma_{22}y_2
\end{cases}
$$
and come to
$$
\begin{pmatrix}
p_1\\ \sigma_{11}y_1+s_{12}y_2  \\p_2\\ \sigma_{21}y_1+\sigma_{22}y_2
\end{pmatrix}
=
\begin{pmatrix}
a_1^*&c_1^*&0&0
\\b_1^*&d_1^*&0&0
\\ 
0&0&a_2^*&c_2^*\\
0&0&b_2^*&d_2^* 
\end{pmatrix}
\begin{pmatrix}
q_1\\ y_1 \\q_2\\ y_2
\end{pmatrix}
$$

{\bf\punct Proof of Theorem \ref{th:main-1}.h.\label{ss:th-h}}
Set $n=2$. We write (\ref{eq:def-2}) as 
\begin{multline*}
\begin{pmatrix}
 q_1\\s_{11}x_1+s_{12}x_2\\q_2\\s_{21}x_1+s_{22}x_2
\end{pmatrix}
=\\=
\begin{pmatrix}
 \lambda_1&0&0&0\\
0&\lambda_1&0&0\\
0&0&\lambda_2&0\\
0&0&0&\lambda_2
\end{pmatrix}^{-1}
\begin{pmatrix}
a_1&b_1&0&0\\
c_1&d_1&0&0\\
0&0&a_2&b_2\\
0&0&c_2&d_2 
\end{pmatrix}
\begin{pmatrix}
 \lambda_1&0&0&0\\
0&\lambda_1&0&0\\
0&0&\lambda_2&0\\
0&0&0&\lambda_2
\end{pmatrix}
\begin{pmatrix}
 p_1\\x_1\\p_2\\x_2
\end{pmatrix}
\end{multline*}

We denote 
$$
y_1=\lambda_1 x_1,\qquad y_2=\lambda_2 x_2,
$$
and write
$$
\begin{pmatrix}
\lambda_1p_1\\
s_{11}y_1 + s_{12}\lambda_1^{-1}\lambda_2 y_2 
\\
\lambda_2p_2 
\\
s_{21}\lambda_1\lambda_2^{-1} y_1 + s_{22} y_2 
\end{pmatrix}
=
\begin{pmatrix}
a_1&b_1&0&0\\
c_1&d_1&0&0\\
0&0&a_2&b_2\\
0&0&c_2&d_2 
\end{pmatrix}
\begin{pmatrix}
\lambda_1 p_1
\\
y_1
\\ \lambda_2p_2\\
y_2
\end{pmatrix}
$$

\section{Language of Grassmannians}

\COUNTERS

Here the rephrase the construction of the previous section.


\smallskip

{\bf\punct Linear relations.} Let $V$, $W$ be linear spaces.
A {\it linear relation} $P:V\tto W$ is a subspace in $V\oplus W$.
The product $QR$ of linear relations $P:V\tto W$, $Q:W\tto Y$
is a linear relation consisting of
$v\oplus y\in V\oplus Y$ such that there exists
$w\in W$ satisfying $v\oplus w\in P$, $w\oplus y\in Q$.


\smallskip

{\sc Example.}  Let $A:V\to W$ is linear operator. Then the graph
$\graph(A)$ of $A$
is a linear relation.

\smallskip


{\bf \punct Eigensurface.}
Denote by $\Gr_{p,q}$ the Grassmannian of $p$-dimensional subspaces
in $\C^{p+q}$. 
Let us keep the notation of the previous section.

Consider the space $\C^n\oplus\C^n$ with coordinates 
$(v_1,\dots,v_n,w_1,\dots, w_n)$. 
 
Fix an $n$-dimensional subspace $L\subset \C^n\oplus\C^n$.
Let 
$$
\sum s_{ij}v_j+\sum \sigma_{ij}w_j=0
$$
be a collection of equations determining $L$.

For definiteness, set $n=3$. We write the equation
\begin{equation}
\begin{pmatrix}
0\\y_1
\\0\\ y_2
 \\0\\ y_3
 \end{pmatrix}
=
\begin{pmatrix}
a_1&b_1&0&0 &0&0\\
c_1&d_1&0&0 &0&0 \\
0&0&a_2&b_2 &0&0 \\
0&0&c_2&d_2 &0&0\\
0&0&0&0&a_3&b_3\\
0&0&0&0&c_3&d_3
\end{pmatrix}
\begin{pmatrix}
0\\ x_1  \\ 0\\ x_2\\0\\x_3
\end{pmatrix}
\label{eq:def-3}
\end{equation}
or, equivalently,
\begin{equation}
\begin{pmatrix}
 y_1\\y_2\\y_3
\end{pmatrix}
=
\begin{pmatrix}
 d_1&0&0\\
0&d_2&0\\
0&0&d_3
\end{pmatrix}
\begin{pmatrix}
 x_1\\x_2\\x_3
\end{pmatrix}
\label{eq:eigen}
\end{equation}

We say that $L\in\Gr_{n,n}$ is contained in the {\it eigensurface}
$\Xi[\frA]$ if there exists a non-zero vector
$(x_1,x_2,x_3,y_1,y_2,y_3)$ satisfying  (\ref{eq:eigen}) and the system
\begin{equation}
\Biggl\{
\sum s_{ij}x_j+\sum \sigma_{ij}y_j=0
\label{eq:L}
\end{equation}
and (\ref{eq:def-3}).


\smallskip

{\bf\punct Characteristic function.} For definiteness, let $n=3$.
Now we write the equation
\begin{equation}
\begin{pmatrix}
 q_1\\y_1\\q_2\\y_2\\q_3\\y_3
\end{pmatrix}
=
\begin{pmatrix}
a_1&b_1&0&0 &0&0\\
c_1&d_1&0&0 &0&0 \\
0&0&a_2&b_2 &0&0 \\
0&0&c_2&d_2 &0&0\\
0&0&0&0&a_3&b_3\\
0&0&0&0&c_3&d_3
\end{pmatrix}
\begin{pmatrix}
 p_1\\x_1\\p_2\\x_2\\p_3\\x_3
\end{pmatrix}
\label{eq:def-25}
\end{equation}
For any $L\in\Gr_{3,3}$
we construct a linear relation
$$
X(\frA;L)=\C^\alpha\oplus\C^\alpha\oplus \C^\alpha
\tto
\C^\alpha\oplus\C^\alpha\oplus \C^\alpha
$$
by the following rule:
$(p_1,p_2,p_3;q_1,q_2,q_3)\in X(\frA;L)$
if there exist $(x_1,x_2,x_3, y_1, y_2, y_)3$ satisfying
(\ref{eq:def-25}) and (\ref{eq:L}).

Note that $X(\frA;L)$ is well-defined also 
for $L$ being in the eigensurface. If 

\begin{theorem}
$$
X(\frA\circ \frP;L)\supset X(\frA;L)\circ X(\frP;L)
.$$
\end{theorem}

{\sc Proof} is the same as in \ref{ss:th-c}.


\smallskip

{\bf\punct Rephrasing of the expansion property.}
We define 
an indefinite Hermitian form in $(\C^\alpha)^n\oplus (\C^\alpha)^n$ by
$$
\cM\bigl(p\oplus q\,\,;\,\,
p'\oplus q'\bigr)=
\sum_{i=1}^n\la p_i,p_i'\ra
-\sum_{i=1}^n\la q_i,q_i'\ra
$$
and a Hermitian form $M$ on $\C^n\oplus \C^n$ 
given by
$$
M(v\oplus w\,;\,
v'\oplus w')=
\sum_{i=1}^n v_i\ov v_i'-\sum_{i=1}^n w_i\ov w_i'
.$$

\begin{theorem}
If $M$ is positive definite on $L$, then $\cM$ is 
negative definite on $X(\frA;L)$. 
\end{theorem}

\section{Conjugacy classes: another example}

\COUNTERS

{\bf \punct  Conjugacy classes.}
Now we consider elements of the group
$\U(\alpha+p+p)$ up to the equivalence
$$
\begin{pmatrix}
a&b_1&b_2
\\
c_1&d_{11}&d_{12}
\\
c_2&d_{21}&d_{22} 
\end{pmatrix}
\sim
\begin{pmatrix}
1&0&0\\
0&u&0\\
0&0&u 
\end{pmatrix}
\begin{pmatrix}
a&b_1&b_2
\\
c_1&d_{11}&d_{12}
\\
c_2&d_{21}&d_{22} 
\end{pmatrix}
\begin{pmatrix}
1&0&0\\
0&u&0\\
0&0&u 
\end{pmatrix}^{-1}
$$
where $u\in\U(p)$.
We denote this set by
$\U(\alpha+2p)//\U(p)$.
We have a well-defined operation
$$
\U(\alpha+2p)\times \U(\alpha+2q)
\to
 \U(\alpha+2p+2q)
$$ 
determined by
\begin{multline*}
\begin{pmatrix}
a&b_1&b_2
\\
c_1&d_{11}&d_{12}
\\
c_2&d_{21}&d_{22} 
\end{pmatrix}
\circ
\begin{pmatrix}
\wt a&\wt b_1&\wt b_2
\\
\wt c_1&\wt d_{11}&\wt d_{12}
\\
\wt c_2&\wt d_{21}&\wt d_{22} 
\end{pmatrix}
:=\\:=
\begin{pmatrix}
a&b_1&0&b_2&0
\\
c_1&d_{11}&0&d_{12}&0
\\
0&0&1&0&0
\\
c_2&d_{21}&0&d_{22}&0
\\
0&0&0&0&1 
\end{pmatrix}
\begin{pmatrix}
\wt a&0&\wt b_1&0&\wt b_2
\\
0&1&0&0&0
\\
\wt c_1&0&\wt d_{11}&0&\wt d_{12}
\\
0&0&0&1&0
\\
\wt c_2&0&\wt d_{21}&0&\wt d_{22} 
\end{pmatrix}
\end{multline*}

We define the characteristic function on the space of
$2\times 2$ matrices. We write the equation
\begin{equation}
\begin{pmatrix}
q\\
s_{11}x_1+s_{12}x_2
\\ 
s_{21}x_1+s_{22}x_2
\end{pmatrix}
=
\begin{pmatrix}
a&b_1&b_2
\\
c_1&d_{11}&d_{12}
\\
c_2&d_{21}&d_{22} 
\end{pmatrix}
\begin{pmatrix}
 p\\x_1\\x_2
\end{pmatrix}
\end{equation}

Then we exclude $x_1$, $x_2$ and get a matrix-valued function on
$\Mat(2)$.

\begin{theorem}
 All the claims of Theorem \ref{th:main-1}
hold except h).
\end{theorem}

Proofs are the same as above.

\section{Example: double cosets}

\COUNTERS

{\bf\punct Double cosets $\U(\infty)\times\dots\times\U(\infty)$
with respect to $\OO(\infty-\alpha)$.}
Now we consider collections of $n$ finite unitary matrices
$$
\Biggl\{\,\begin{pmatrix}a_j&b_j\\c_j&d_j\end{pmatrix}\,
\Biggr\}
$$
of the size $(\alpha+\infty)\times(\alpha+\infty)$
determined up to the equivalence
\begin{multline}
\Biggl\{\,\begin{pmatrix}a_1&b_1\\c_1&d_1\end{pmatrix},
\dots,
\begin{pmatrix}a_n&b_n\\c_n&d_n\end{pmatrix}
\Biggr\}
\sim\\ \sim
\Biggl\{\,
\begin{pmatrix} 1&0\\0&u \end{pmatrix}
\begin{pmatrix}a_1&b_1\\c_1&d_1\end{pmatrix}
\begin{pmatrix} 1&0\\0&v \end{pmatrix}
,\dots,
\begin{pmatrix} 1&0\\0&u \end{pmatrix}
\begin{pmatrix}a_n&b_n\\c_n&d_n\end{pmatrix}
\begin{pmatrix} 1&0\\0&v \end{pmatrix}\,
\Biggr\}
\label{eq:class-last}
\end{multline}
where $u$, $v$ are finite real orthogonal matrices,
$u=u^{t-1}$, $v=v^{t-1}$.
The multiplication of colligations is defined as above.

\smallskip


{\bf \punct Characteristic functions.}
For definiteness, set $n=2$.
We write the equation
\begin{equation}
\begin{pmatrix}
 q_1^+\\y_1^+\\q_2^+\\y_2^+\\ q_1^-&\\y_1^-&\\q_2^-\\ y_2^-
\end{pmatrix}
=
\begin{pmatrix}
\begin{pmatrix}
a_1&b_1\\c_1&d_1
\end{pmatrix} &&&
\\
&\begin{pmatrix}
a_2&b_2\\c_2&d_2
\end{pmatrix}&&
\\
&&\begin{pmatrix}
a_1&b_1\\c_1&d_1
\end{pmatrix}^{t-1}&
\\
&&&\begin{pmatrix}
a_2&b_2\\c_2&d_2
\end{pmatrix}^{t-1}
\end{pmatrix}
\begin{pmatrix}
 p_1^+\\x_1^+\\p_2^+\\x_2^+\\ p_1^-&\\x_1^-&\\p_2^-\\ x_2^-
\end{pmatrix}
\label{eq:def-last}
\end{equation}

Next, take two $n\times n$ matrices $S$, $R$ (in our case $n=2$)
 and assume that
$x$, $y$ satisfy the equations
\begin{equation}
\begin{cases}
y_1^+=s_{11}y_1^-+s_{12}y_2^-
\\
y_2^+=s_{21}y_1^-+s_{22}y_2^-
\end{cases}
\qquad
\begin{cases}
x_1^-=r_{11}x_1^++r_{12}x_2^+
\\
x_2^-=r_{21}x_1^++r_{22}x_2^+
\end{cases}
\label{eq:rs}
\end{equation}

Now we exclude variables $y_1^+$, $y_2^+$, $x_1^-$, $x_2^-$
and come to a dependence of the form
$$
\begin{pmatrix}
 q_1^+\\q_2^+\\q_1^-\\q_2^-
\end{pmatrix}
=
\chi(\frA;S,R)
\begin{pmatrix}
 p_1^+\\p_2^+\\p_1^-\\p_2^-
\end{pmatrix}
,
$$
where $\chi(\frA;T,R)$ is a function of two variables
 $T$, $R\in\Mat(n)$ taking values in
$\Mat(2n\alpha)$.

\begin{theorem}
 \label{th:main-last}
a)  $\chi(\frA;S,R)$  is a meromorphic matrix-valued function.

\smallskip

b) $\chi(\frA;S,R)$ depends only of the equivalence class
(\ref{eq:class-last}) and not on matrices themselves.

\smallskip

c)
$\chi(\frA\circ\frP;S,R)=\chi(\frA;S,R)\,\chi(\frA;S,R)$
\end{theorem} 

Proofs are similar to proofs of corresponding statements
of Theorem \ref{th:main-1}.


\smallskip

Theorem \ref{th:main-1} also can be extended in a straightforward way.
For definiteness, set $n=2$.

\begin{theorem} Let $\lambda_1$, $\lambda_2\in\C^*$.
\begin{multline*}
\begin{pmatrix}
 \lambda_1&0&0&0\\
0&\lambda_2&0&0\\
0&0&\lambda_1^{-1}&0\\
0&0&0&\lambda_2^{-1}
\end{pmatrix}
\chi(\frA;S,R)
\begin{pmatrix}
 \lambda_1&0&0&0\\
0&\lambda_2&0&0\\
0&0&\lambda_1^{-1}&0\\
0&0&0&\lambda_2^{-1}
\end{pmatrix}^{-1}
=\\=
\chi\left[\frA; 
\begin{pmatrix}
 \lambda_1&0\\
0&\lambda_2
\end{pmatrix}
S
\begin{pmatrix}
 \lambda_1&0\\
0&\lambda_2
\end{pmatrix}^{-1}\,\,
,\,\, 
\begin{pmatrix}\lambda_1&0\\
0&\lambda_2
\end{pmatrix}
 R 
\begin{pmatrix}
 \lambda_1&0\\
0&\lambda_2
\end{pmatrix}^{-1}\right]
\end{multline*}

\end{theorem}


{\bf\punct The analog of expansion property.}
We define an Hermitian form $\cM$
in the space $V:=(\C^\alpha)^n\oplus (\C^\alpha)^n$
by 
$$
\cM\bigl(p_1^+,\dots,p_n^+,p_1^-,\dots, p_n^-\,;\,
\wt p_1^+,\dots, \wt p_n^+,\wt p_1^-,\dots,\wt p_n^-)=
\sum_{j=1}^n \la p_j^+, \wt p_j^+\ra
-
\sum_{j=1}^n \la p_j^- ,\wt p_j^-\ra
.
$$

\begin{theorem}
a) Let $\|S\|<1$, $\|R\|<1$. Let $q=\chi(\frA;S,R)p$.
Then
$$M(q,q)>M(p,p).$$ 

b) Let $S$, $R$ be unitary, let $(S,R)$ be a point of holomorpy 
of $\chi(\frA;S,R)$. Then $\chi(\frA;S,R)$ is contained in the
pseudo-unitary group
$\U(n\alpha,n\alpha)$, i.e.,
 $$M(q,q)=M(p,p) \qquad\text{if $q=\chi(\frA;S,R)p$}
.$$ 

c) The following Riemann--Scwartz type identity holds
$$
\chi(\frA;S^{\square-1},R^{\square-1})=
\chi(\frA;S,R)^{\square-1}
$$
where $^\square$ denotes  the adjoint operator with respect 
to the form $M$.
\end{theorem}

{\sc Proof.} We prove a). We introduce the following
 indefinite Hermitian form $\bfM$ on the space
$V\oplus \ell_2$,
$$
\bfM(p^+,p^-,x^+,x^-\, ;\,\wt p^+,\wt p^-,\wt x^+,\wt x^-)
=\la p^+,\wt p^+\ra+\la x^+,\wt x^+\ra
- \la p^-,\wt p^-\ra - \la x^-,\wt x^-\ra
.$$

The big matrix in (\ref{eq:def-last})
preserves this form. Therefore, for $(p,x)$ and $(q,y)$
related by (\ref{eq:def-last}) the following identity
holds
$$
\|p^+\|^2-\|p^-\|^2+\|x^+\|^2-\|x^-\|^2
=\|q^+\|^2-\|q^-\|^2+\|y^+\|^2-\|y^-\|^2
.$$
Now we note that $x_-=Rx_+$, $y+=Sy_-$. Therefore
$$
\|x^-\|^2<\|x^+\|^2 \qquad \|y^+\|^2<\|y^-\|^2
$$
This implies the required statement.

\smallskip

{\sc Remark.} In Theorem \ref{th:main-1}.e
 $\chi^{-1}$ is  holomorphic in the matrix ball,
this is not valid here. \hfill $\square$

\smallskip


{\bf \punct An additional symmetry.}
In our case, there is an additional property of characteristic functions
comparing(?) Theorem \ref{th:main-1}.

Let us introduce the following skew-symmetric bilinear form
$\Lambda$ in $V$,
$$
\Lambda(p,\wt p)=\sum_{j=1}^n 
\bigl( (p_j^+ ,\wt p_j^-) - (p_j^- \wt p_j^+)\bigr).
.$$
Denote by $^\triangle$ the transposition with respect to the form
$\Lambda$,
$$
\Lambda(g p,\wt p)= \Lambda(p,g^\triangle \wt p)
.$$

\begin{theorem} a)
$$ 
\chi(\frA;S^t, R^t)=\chi(\frA;S,R)^{\triangle-1}
.$$

b) In particular, if matrices $S$, $R$ are symmetric, then 
$\chi(\frA;T,S)$ is contained in the symplectic group
 $\Sp(2n\alpha,\C)$.
\end{theorem}

{\sc Proof.} We introduce a skew-symmetric bilinear form 
$\cL$ in $V\oplus \ell_2$ by
$$
\cL(p\oplus x,\wt p\oplus \wt x)=
\sum_{j=1}^n \bigl((p_j^+, \wt p_j^-) - (p_j^- ,\wt p_j^+)\bigr)+
\sum_{j=1}^n \bigl((x_j^+, \wt x_j^-) - (x_j^- ,\wt x_j^+)\bigr)
.$$
Let $p\oplus x$, $q\oplus y$ 
(and also $\wt p\oplus \wt x$, $\wt q\oplus \wt y$)
 satisfy (\ref{eq:def-last}). Let $x$, $y$ be connected by (\ref{eq:rs}) and
$\wt x$, $\wt y$ be connected by the same relation (\ref{eq:rs}), where $R$, $S$
are replaced by $R^t$, $S^t$.
Then
$$
\cL(p\oplus x,\wt p\oplus \wt x)=\cL(q\oplus y,\wt q\oplus \wt y)
$$
Let 
$$
y^+=Sy^-, \quad \wt y^+=S^t y^-.
\quad 
x^-=R x^+,\quad \wt x^-=R^t x^+
.$$
Then
\begin{align*}
(y^+,\wt y^-)-(y^-,\wt y^+)=(Sy^-,\wt y^-)-(y^-, S^t \wt y^+)=0;
\\
(x^+,\wt x^-)-(x^-,\wt x^+)=(x^+,R\wt x^+)-(R^tx^+,\wt x^+)=0
.\end{align*}
Therefore
$$
\Lambda(q,\wt q)=\Lambda(p,\wt p).
$$

{\tt Math.Dept., University of Vienna,

 Nordbergstrasse, 15,
Vienna, Austria

\&

Institute for Theoretical and Experimental Physics,

Bolshaya Cheremushkinskaya, 25, Moscow 117259,
Russia

\&

Mech.Math. Dept., Moscow State University,
Vorob'evy Gory, Moscow

e-mail: neretin(at) mccme.ru

URL:www.mat.univie.ac.at/$\sim$neretin

wwwth.itep.ru/$\sim$neretin

\end{document}